\documentclass[12pt,a4paper]{article}
\usepackage{amsmath,amsthm,amssymb}
\usepackage[truedimen,margin=25truemm]{geometry}
\usepackage{mathrsfs}
\usepackage{bm}
\usepackage[utf8]{inputenc}
\usepackage{mathtools}
\usepackage[hidelinks]{hyperref}

\theoremstyle{definition}
\newtheorem{definition}{Definition}
\newtheorem{thm}[definition]{Theorem}

\newtheorem{lem}[definition]{Lemma}
\newtheorem{cor}[definition]{Corollary}

\numberwithin{definition}{section}

\newcommand{\diag}{{\rm diag}}

\newcommand{\Mat}{{\rm Mat}}

\DeclareMathOperator{\SL}{SL}

\DeclareMathOperator{\sgn}{sgn}

\DeclareMathOperator{\Tr}{Tr}

\DeclareMathOperator{\Ker}{Ker}

\DeclareMathOperator{\Core}{Core}
\DeclareMathOperator{\Cdet}{cdet}

\hypersetup{
  pdftitle={Invariant rings of Frobenius polarization operators and normal cores of subgroups},
  pdfauthor={Naoya Yamaguchi}
}

\begin{document}

\title{Invariant rings of Frobenius polarization operators and normal cores of subgroups}
\author{Naoya Yamaguchi}
\date{}
\maketitle

\begin{abstract}
For a proper subgroup $H$ of a finite group $G$, we determine the common kernel,
or invariant ring, of the Frobenius polarization operators corresponding to
the elements of $G\setminus H$.
We show that the invariant rings defined by the left and right polarization
operators coincide and depend only on the normal core $N=\Core_G(H)$.
We construct polynomials $P_\Omega$ indexed by the orbits of the conjugation
action of $G$ on $\widehat N$ and prove that they form an algebraically
independent generating set for the invariant ring.
We also show that inclusion of these invariant rings as subrings of
$\mathbb C[\bm x]$ is equivalent to inclusion of the corresponding normal cores.
For an arbitrary subgroup $H$, we give a factorization of the group determinant
in terms of these polynomials. This factorization generalizes Frobenius'
group determinant theorem and reduces to its factorization when $H=G$.
When $G$ is abelian, each $P_\Omega$ is a specialization of the group determinant
of the quotient group $G/H$.
Finally, we establish a Capelli identity for each $P_\Omega$ and compute the
coefficient polynomial obtained by applying a constant-coefficient differential
operator to its powers.
\end{abstract}

\section{Introduction}
In this paper, we describe the invariant rings of Frobenius polarization
operators associated with subgroups of finite groups.

\subsection*{Main results}
Let $G$ be a finite group, and let $\widehat G$ be a complete set of
representatives of the equivalence classes of finite-dimensional irreducible
unitary representations of $G$ over $\mathbb C$.
Let $\bm x\coloneqq\{x_g\}_{g\in G}$ be a family of independent commuting
variables indexed by the elements of $G$.
For each $g\in G$, define the Frobenius polarization operators by
\[
D_g\coloneqq\sum_{h\in G}x_{g^{-1}h}\frac{\partial}{\partial x_h},
\qquad
\Delta_g\coloneqq\sum_{h\in G}x_{hg^{-1}}\frac{\partial}{\partial x_h}.
\]
Both operators are derivations of $\mathbb C[\bm x]$.
Their kernels are denoted by
\[
\begin{aligned}
\Ker D_g&\coloneqq\{f(\bm x)\in\mathbb C[\bm x]\mid D_gf(\bm x)=0\},\\
\Ker\Delta_g&\coloneqq\{f(\bm x)\in\mathbb C[\bm x]\mid\Delta_gf(\bm x)=0\}.
\end{aligned}
\]
For a subgroup $H$ of $G$, define the common kernels, or invariant rings, by
\[
K_H^D\coloneqq\bigcap_{g\in G\setminus H}\Ker D_g,
\qquad
K_H^\Delta\coloneqq\bigcap_{g\in G\setminus H}\Ker\Delta_g.
\]
When $H=G$, we use the convention that the common kernel of an empty family is
the entire polynomial ring, so that
\[
K_H^D=K_H^\Delta=\mathbb C[\bm x].
\]
These invariant rings are subrings of $\mathbb C[\bm x]$.

For a proper subgroup $H$, we establish a first fundamental theorem,
which describes generators, and a second fundamental theorem,
which describes the relations among these generators.
We also show that the invariant ring, regarded as a subring of
$\mathbb C[\bm x]$, uniquely determines the normal core of $H$.
The normal core of $H$ in $G$ is defined by
\[
N\coloneqq\Core_G(H)=\bigcap_{t\in G}t^{-1}Ht.
\]
It is the largest normal subgroup of $G$ contained in $H$.
Similarly, let $\widehat N$ be a complete set of representatives of the
equivalence classes of finite-dimensional irreducible unitary representations
of $N$ over $\mathbb C$.
Identifying a representation with its equivalence class, we define the
conjugation action of $G$ on $\widehat N$ by
\[
(g\cdot\sigma)(n)\coloneqq\sigma(g^{-1}ng)
\qquad(g\in G,\ n\in N).
\]
Let $\mathcal O_N$ denote the set of orbits of this action.
The following theorem gives generators for the invariant rings
$K_H^D$ and $K_H^\Delta$.

\begin{thm}[First fundamental theorem]\label{thm:1.1}
Let $H$ be a proper subgroup of a finite group $G$, and set
$N\coloneqq\Core_G(H)$. Then
\[
K_H^D=K_H^\Delta=K_N^D=K_N^\Delta
=\mathbb C[P_\Omega\mid\Omega\in\mathcal O_N].
\]
\end{thm}

The polynomials $P_\Omega$ are defined immediately after
Theorem~\ref{thm:1.5}.
By Theorem~\ref{thm:1.1} and the definition for $H=G$, the invariant rings
$K_H^D$ and $K_H^\Delta$ coincide. Henceforth, we write
$K_H\coloneqq K_H^D=K_H^\Delta$.
The next theorem asserts that the generators
$\{P_\Omega\}_{\Omega\in\mathcal O_N}$ in Theorem~\ref{thm:1.1}
are algebraically independent.

\begin{thm}[Second fundamental theorem]\label{thm:1.2}
Let $H$ be a subgroup of a finite group $G$, and set
$N\coloneqq\Core_G(H)$. Then the polynomials
$\{P_\Omega\}_{\Omega\in\mathcal O_N}$ are algebraically independent
over $\mathbb C$.
\end{thm}

Thus, there are no nontrivial polynomial relations among these polynomials.
It is in this sense that we refer to Theorem~\ref{thm:1.2} as the second
fundamental theorem.

When $H$ is a proper subgroup of $G$, Theorems~\ref{thm:1.1}
and~\ref{thm:1.2} show that $K_H$ is isomorphic to a polynomial ring in
$|\mathcal O_N|$ variables over $\mathbb C$.
A consequence of Brauer's permutation lemma~\cite[Corollary~6.33]{MR460423}
identifies $|\mathcal O_N|$ with the number of conjugacy classes of $G$
contained in $N$.
A conjugacy class of $G$ is contained in $H$ if and only if it is contained
in $N$. We therefore obtain the following corollary.

\begin{cor}\label{cor:1.3}
Let $H$ be a proper subgroup of a finite group $G$, and set
$N\coloneqq\Core_G(H)$. The number of algebraically independent generators
of the invariant ring $K_H$ over $\mathbb C$ is
\[
|\mathcal O_N|
=|\{C\mid C\text{ is a conjugacy class of }G,\ C\subset H\}|.
\]
\end{cor}

The following theorem asserts that $K_H$, regarded as a subring of
$\mathbb C[\bm x]$, uniquely determines the normal core $\Core_G(H)$.

\begin{thm}\label{thm:1.4}
Let $H_1,H_2$ be proper subgroups of a finite group $G$, and set
\[
N_i\coloneqq\Core_G(H_i),
\qquad
K_{H_i}\coloneqq K_{H_i}^D=K_{H_i}^\Delta
\qquad(i=1,2).
\]
Then
\[
K_{H_1}\subset K_{H_2}
\iff
N_1\subset N_2.
\]
In particular,
\[
K_{H_1}=K_{H_2}
\iff
\Core_G(H_1)=\Core_G(H_2).
\]
\end{thm}

We next factor the group determinant
\[
\Theta_G(\bm x)\coloneqq\det(x_{gh^{-1}})_{g,h\in G}
\]
as a product of suitable powers of the polynomials $P_\Omega$.
For each orbit $\Omega\in\mathcal O_N$, set
\[
I_\Omega\coloneqq
\{\varphi\in\widehat G\mid
\varphi|_N\text{ has an irreducible constituent in }\Omega\}.
\]
The family $\{I_\Omega\}_{\Omega\in\mathcal O_N}$ is a partition of
$\widehat G$, and every $I_\Omega$ is nonempty
(Lemma~\ref{lem:2.20}). Define
\[
\delta_\Omega\coloneqq\gcd\{d_\varphi\mid\varphi\in I_\Omega\},
\qquad d_\varphi\coloneqq\deg\varphi.
\]

\begin{thm}[A generalization of Frobenius' group determinant theorem]\label{thm:1.5}
Let $H$ be a subgroup of a finite group $G$, and set $N=\Core_G(H)$.
For the polynomials $\{P_\Omega\}_{\Omega\in\mathcal O_N}$ defined below,
we have
\[
\Theta_G(\bm x)=\prod_{\Omega\in\mathcal O_N}P_\Omega^{\delta_\Omega}.
\]
\end{thm}

We now define the polynomials $P_\Omega$.
When $H$ is a proper subgroup, these are the generators in
Theorem~\ref{thm:1.1}.
For each $\varphi\in\widehat G$, set
\[
X_\varphi\coloneqq\sum_{g\in G}\varphi(g)x_g.
\]
For each $\Omega\in\mathcal O_N$, define
\[
P_\Omega\coloneqq
\prod_{\varphi\in I_\Omega}(\det X_\varphi)^{k_\varphi},
\qquad
k_\varphi\coloneqq\frac{d_\varphi}{\delta_\Omega}
\quad(\varphi\in I_\Omega).
\]
These definitions also apply when $H=G$.

Let $e$ denote the identity element of $G$.
First suppose that $G\ne\{e\}$, and apply Theorem~\ref{thm:1.1}
to the proper subgroup $\{e\}$.
In this case, $N=\{e\}$, and $\mathcal O_{\{e\}}$ consists of a single
orbit $\Omega_0$.
We have $I_{\Omega_0}=\widehat G$. Since the trivial representation of $G$
has degree $1$,
\[
\delta_{\Omega_0}
=\gcd\{d_\varphi\mid\varphi\in\widehat G\}
=1.
\]
Thus, Theorems~\ref{thm:1.1} and~\ref{thm:1.5} give
\[
K_{\{e\}}
=\mathbb C[P_{\Omega_0}]
=\mathbb C[\Theta_G(\bm x)].
\]
Applying Theorem~\ref{thm:1.4} to a proper subgroup $H$ and to $\{e\}$
proves Corollary~\ref{cor:1.6} below for proper subgroups.
We also check the case $H=G$.
If $G=\{e\}$, then $\Theta_G(\bm x)=x_e$, and both conditions hold.
If $|G|\ge2$, then $\Theta_G$ is homogeneous of degree $|G|$, so that
$x_e\notin\mathbb C[\Theta_G]$. Hence
$K_G=\mathbb C[\bm x]\ne\mathbb C[\Theta_G]$.
Moreover, $\Core_G(G)=G\ne\{e\}$.
Consequently, the following statement also holds when $H=G$.

\begin{cor}\label{cor:1.6}
Let $H$ be a subgroup of a finite group $G$. Then
$K_H=\mathbb C[\Theta_G(\bm x)]$ if and only if $\Core_G(H)=\{e\}$.
\end{cor}

By definition, each polynomial $P_\Omega$ has the determinantal expression
\[
P_\Omega
=\det\left(\sum_{g\in G}x_g\rho_\Omega(g)\right),
\qquad
\rho_\Omega\coloneqq
\bigoplus_{\varphi\in I_\Omega}\varphi^{\oplus k_\varphi}.
\]

The next theorem shows that, when $G$ is abelian, each $P_\Omega$ is a
specialization of the group determinant of the quotient group $G/H$.

\begin{thm}\label{thm:1.7}
Let $H$ be a subgroup of a finite abelian group $G$.
For each $\sigma\in\widehat H$, choose $\varphi_\sigma\in\widehat G$
such that $\varphi_\sigma|_H=\sigma$.
For each $tH\in G/H$, set
\[
y_{tH}^\sigma\coloneqq
\sum_{g\in tH}\varphi_\sigma(g)x_g
=\sum_{h\in H}\varphi_\sigma(th)x_{th},
\]
and write $\bm y^\sigma\coloneqq\{y_{tH}^\sigma\}_{tH\in G/H}$.
Then, for every $\sigma\in\widehat H$,
\[
P_{\{\sigma\}}=\Theta_{G/H}(\bm y^\sigma).
\]
\end{thm}

When $H$ is a proper subgroup of $G$, the next corollary follows from
Theorems~\ref{thm:1.1}, \ref{thm:1.2}, \ref{thm:1.5}, and~\ref{thm:1.7}.
When $H=G$, we have
\[
\Theta_{G/H}(\bm y^\sigma)
=\sum_{g\in G}\sigma(g)x_g
\qquad(\sigma\in\widehat G).
\]
By character orthogonality, these linear forms constitute a linear coordinate
system for $\mathbb C[\bm x]$. Thus, the same conclusion holds in this case.

\begin{cor}\label{cor:1.8}
Let $H$ be a subgroup of a finite abelian group $G$. Then
\[
K_H^D=K_H^\Delta
=\mathbb C\bigl[
\Theta_{G/H}(\bm y^\sigma)
\mid\sigma\in\widehat H
\bigr].
\]
These $|H|$ generators are algebraically independent over $\mathbb C$,
and each has degree $[G:H]$.
Moreover,
\[
\Theta_G(\bm x)
=\prod_{\sigma\in\widehat H}\Theta_{G/H}(\bm y^\sigma).
\]
\end{cor}

Finally, we give a Capelli identity for the invariants $P_\Omega$.
For each $\varphi\in\widehat G$, set
\[
\partial_\varphi\coloneqq
\sum_{g\in G}\overline{\varphi(g)}\frac{\partial}{\partial x_g},
\quad
\alpha_\varphi\coloneqq\frac{|G|}{d_\varphi},
\quad
\natural_\varphi\coloneqq\diag(d_\varphi-1,d_\varphi-2,\ldots,0).
\]
For an indeterminate $u$ commuting with all operators, define
\[
\begin{aligned}
C_\varphi(u)
&\coloneqq\Cdet\left(
{}^tX_\varphi\partial_\varphi
+\alpha_\varphi(uI_{d_\varphi}+\natural_\varphi)\right)\\
&=\Cdet\left(
\sum_{g\in G}\overline{\varphi(g)}\Delta_g
+\alpha_\varphi(uI_{d_\varphi}+\natural_\varphi)\right).
\end{aligned}
\]
Here, $\Cdet$ denotes the column determinant defined in
Section~\ref{sec:2.2}, and $I_{d_\varphi}$ is the identity matrix of size
$d_\varphi$.
The equality of the two expressions follows from Lemma~\ref{lem:2.9}.
For each $k\in\mathbb Z_{\ge0}$, set
\[
C_\varphi^{(k)}(u)\coloneqq
\begin{cases}
1,&k=0,\\
C_\varphi(u)C_\varphi(u-1)\cdots C_\varphi(u-k+1),&k\ge1,
\end{cases}
\]
where $1$ denotes the identity operator.

\begin{thm}\label{thm:1.9}
Let $H$ be a subgroup of a finite group $G$, and set $N=\Core_G(H)$.
For each orbit $\Omega\in\mathcal O_N$, define the constant-coefficient
differential operator
\[
\mathcal D_\Omega\coloneqq
\prod_{\varphi\in I_\Omega}(\det\partial_\varphi)^{k_\varphi}.
\]
Then, as differential operators on $\mathbb C[\bm x]$,
\[
P_\Omega\mathcal D_\Omega
=\prod_{\varphi\in I_\Omega}C_\varphi^{(k_\varphi)}(0).
\]
Moreover, for every $s\in\mathbb Z_{\ge0}$,
\[
\mathcal D_\Omega(P_\Omega^{s+1})=b_\Omega(s)P_\Omega^s,
\]
where
\[
b_\Omega(s)\coloneqq
\prod_{\varphi\in I_\Omega}
\left\{
\alpha_\varphi^{d_\varphi k_\varphi}
\prod_{a=1}^{k_\varphi}\prod_{i=0}^{d_\varphi-1}
(k_\varphi s+a+i)
\right\}.
\]
\end{thm}
\subsection*{Relation to previous work}
We now discuss the relation between our main results and previous work.

\paragraph{Frobenius' group determinant theorem and polarization operators.}
Frobenius~\cite{Frobenius1968gruppen} introduced the polarization operators
$D_g$ and $\Delta_g$ in his study of the irreducible factorization
\[
\Theta_G(\bm x)
=\prod_{\varphi\in\widehat G}(\det X_\varphi)^{d_\varphi}
\]
of the group determinant over $\mathbb C$.
Let $\chi_\varphi$ denote the character of $\varphi$.
The relation between these operators and the irreducible factors is expressed by
\[
D_g\det X_\varphi
=\Delta_g\det X_\varphi
=\chi_\varphi(g)\det X_\varphi
\]
(see, for example, \cite[p.~226]{MR803326} and
\cite[Sections~1 and~5]{Yamaguchi13022026}).
By this formula and the Leibniz rule, a product of nonnegative integral powers
of the irreducible factors is an eigenvector of every polarization operator.
Its eigenvalue function is the corresponding nonnegative integral linear
combination of irreducible characters.

The factorization in Theorem~\ref{thm:1.5} follows by grouping the factors
in Frobenius' irreducible factorization according to the sets $I_\Omega$
and factoring out the common divisor $\delta_\Omega$ of the exponents.
For a proper subgroup $H$, we prove that the polynomials $P_\Omega$ obtained
by this grouping and normalization of the exponents form an algebraically
independent generating set for the common kernel $K_H$ in the entire
polynomial ring $\mathbb C[\bm x]$.
When $H=G$, the normal core is $N=G$, and the conjugation action on
$\widehat G$ is trivial. Hence
\[
I_{\{\varphi\}}=\{\varphi\},\qquad
\delta_{\{\varphi\}}=d_\varphi,\qquad
P_{\{\varphi\}}=\det X_\varphi,
\]
and the identity in Theorem~\ref{thm:1.5} reduces to Frobenius'
group determinant formula.

\paragraph{Kernels of individual polarization operators.}
For diagonal derivations, which send each variable to a scalar multiple of
itself, Nowicki and Strelcyn~\cite[Theorem~1]{NowickiStrelcyn1995} proved that,
over a field $K$ of characteristic zero, the minimum number of generators
of the kernel as a $K$-algebra can be any nonnegative integer, even when
the number of variables is fixed at any value at least three.
The polarization operator corresponding to a generator of a cyclic group
has also been studied as a cyclotomic derivation.
For a cyclic group of order $n\ge3$, Moulin Ollagnier and Nowicki
\cite[Propositions~2.3 and~3.1]{MoulinOllagnierNowicki2013} described a
generating set for its kernel consisting of monomials in Fourier-transformed
variables. Their description uses the minimal nonzero solutions of
nonnegative integral vanishing relations among the $n$th roots of unity,
and they proved that the number of generators is minimal.
They also showed in \cite[Theorem~3.7]{MoulinOllagnierNowicki2013} that this
kernel is isomorphic to a polynomial ring over $\mathbb C$ if and only if
$n$ is a prime power.

Combining Theorems~1 and~2 of \cite{YamaguchiOchiaiYamaguchi2025} gives,
when the order $n\ge2$ of the cyclic group has at most two distinct prime
divisors, generators for the kernel as specializations of circulant
determinants of prime order, together with the ideal of relations among
these generators.
Theorem~3 of \cite{YamaguchiOchiaiYamaguchi2025} shows that, when $n$ has
at least three distinct prime divisors, the specializations of circulant
determinants of prime order constructed there do not generate the entire kernel.
For the kernel of the polarization operator corresponding to each element
of an arbitrary finite group, \cite[Theorems~1 and~3]{Yamaguchi2026Kernels}
describes a generating set of minimum cardinality using the discrete
Fourier transform on each coset of the cyclic subgroup generated by that
element and nonnegative integral vanishing relations among roots of unity.

These results concern the generators and ring structure of kernels of
individual derivations.
In the present paper, for an arbitrary finite group $G$ and a proper
subgroup $H$, we study the common kernel of the polarization operators
corresponding to all elements of $G\setminus H$.
Theorems~\ref{thm:1.1} and~\ref{thm:1.2} show that this common kernel is
isomorphic to a polynomial ring over $\mathbb C$ and provide explicit
algebraically independent generators.

\paragraph{Results in the subring generated by the irreducible factors and their extension.}
Let
\[
R\coloneqq\mathbb C[\det X_\varphi\mid\varphi\in\widehat G]
\]
be the subring generated by the irreducible factors of the group determinant.
Within this subring, \cite[Theorem~2]{Yamaguchi13022026} describes the common
kernel of the polarization operators corresponding to representatives of
all conjugacy classes other than the identity class.
Since characters are class functions, polarization operators corresponding
to conjugate elements have the same action on $R$.
Thus, in our notation, that result is
\[
R\cap K_{\{e\}}^D
=R\cap K_{\{e\}}^\Delta
=\mathbb C[\Theta_G(\bm x)].
\]
When $G$ is abelian, character orthogonality gives $R=\mathbb C[\bm x]$,
so this also determines the common kernel in the entire polynomial ring
\cite[Corollary~3]{Yamaguchi13022026}.
Theorem~6 of \cite{Yamaguchi13022026} also describes the kernels of the
individual polarization operators for the symmetric group $S_3$ within $R$.

Furthermore, \cite[Theorem~1]{YI2026} uses a Capelli identity of group
determinant type to prove, for every finite group, the equality
\[
K_{\{e\}}^D
=K_{\{e\}}^\Delta
=\mathbb C[\Theta_G(\bm x)]
\]
in the entire polynomial ring $\mathbb C[\bm x]$.
We extend this description of the common kernel to arbitrary proper
subgroups $H$.
Theorem~\ref{thm:1.1} shows that the common kernel depends only on the
normal core $N=\Core_G(H)$ and gives an explicit generating set.
Theorem~\ref{thm:1.2} proves that these generators are algebraically
independent, and Theorem~\ref{thm:1.4} shows that inclusion of the common
kernels inside the same polynomial ring $\mathbb C[\bm x]$ is equivalent
to inclusion of the normal cores.
Thus, it is the structure of the common kernel as a subring, including its
embedding in $\mathbb C[\bm x]$, that uniquely determines the normal core.
Corollary~\ref{cor:1.6} further shows that the common kernel is
$\mathbb C[\Theta_G(\bm x)]$ if and only if the normal core is trivial.

\paragraph{Relation to classical invariant theory.}
For the action of the special linear group $\SL_d(\mathbb C)$ on the
matrix space $\Mat_d(\mathbb C)$ by left multiplication, the invariant ring is
\[
\mathbb C[\Mat_d(\mathbb C)]^{\SL_d(\mathbb C)}
=\mathbb C[\det]
\]
\cite[Section~1.2, Example~2]{KraftProcesi1996}.
By transposition, the same equality holds for the action by right multiplication.
Lemmas~\ref{lem:2.7} and~\ref{lem:2.8} show that every nonzero common
eigenvector of the family $\{\Delta_g\mid g\in G\}$ is a scalar multiple
of a product of nonnegative integral powers of the determinants
$\det X_\varphi$.
This description can also be obtained by applying the classical description
of the invariant ring above to each matrix block and then imposing the
eigenvector conditions for the Euler operators.
Indeed, a common eigenvector is annihilated by the infinitesimal action of
the special linear group on each matrix block, and the eigenvalues of the
diagonal Euler operators determine the exponents of the determinants.
Here, we prove this description of the common eigenspaces directly using
polarization operators and the Capelli identity.
We then determine the exponent vectors of the products of determinants
occurring in the common kernel by requiring the eigenvalue function to
vanish on $G\setminus N$.

\paragraph{Clifford theory and the exponents of the generators.}
The partition of the irreducible representations into the sets
$\{I_\Omega\}_{\Omega\in\mathcal O_N}$ can be explained by Clifford
theory~\cite{Clifford1937}, which describes restrictions to normal subgroups.
By Clifford's theorem~\cite[Theorem~6.2]{MR460423}, for each
$\varphi\in I_\Omega$ there is a positive integer $e_\varphi$ such that
\[
\varphi|_N\simeq
\bigoplus_{\sigma\in\Omega}\sigma^{\oplus e_\varphi}.
\]
All irreducible representations in an orbit $\Omega$ have the same degree.
Hence, for any $\sigma\in\Omega$,
\[
d_\varphi=e_\varphi|\Omega|d_\sigma.
\]
The definitions of the exponents therefore give
\[
\begin{aligned}
\delta_\Omega
&=|\Omega|d_\sigma\gcd\{e_\psi\mid\psi\in I_\Omega\},\\
k_\varphi
&=\frac{e_\varphi}{\gcd\{e_\psi\mid\psi\in I_\Omega\}}.
\end{aligned}
\]
Thus, the family of exponents $\{k_\varphi\}_{\varphi\in I_\Omega}$ is
obtained by dividing the restriction multiplicities
$\{e_\varphi\}_{\varphi\in I_\Omega}$ by their greatest common divisor.
We use these orbits and multiplicities from restriction theory to describe
algebraically independent generators of the common kernel associated with
a proper subgroup.

\paragraph{Factorizations of the group determinant using subgroups.}
A generalization expressing the group determinant of a finite abelian
group as a specialization of the group determinant of a subgroup was given
in \cite[Theorem~1.4]{Yamaguchi2017Dedekind}.
Theorem~1.4 of \cite{Yamaguchi2019Study} gives a Dedekind-type factorization
over the group algebra for an arbitrary finite group and an abelian subgroup.
For finite abelian groups, the generalization of Dedekind's theorem in
\cite[Theorem~1.1]{MR4814714} gives, in our notation, the factorization
\[
\Theta_G(\bm x)
=\prod_{\sigma\in\widehat H}\Theta_{G/H}(\bm y^\sigma).
\]
This is the last identity in Corollary~\ref{cor:1.8}.
For each factor in this factorization, \cite[Theorem~5]{Yamaguchi13022026}
shows that
\[
D_g\Theta_{G/H}(\bm y^\sigma)=0
\iff g\in G\setminus H.
\]
Theorem~\ref{thm:1.7} identifies these factors with the generators
$P_{\{\sigma\}}$.
Moreover, Corollary~\ref{cor:1.8} shows that they generate the entire
invariant ring $K_H$ and are algebraically independent.
Thus, the factors in this known group determinant factorization are
identified as generators of the common kernel associated with the subgroup.

\paragraph{Relation to Capelli identities.}
Umeda~\cite{MR3837932} studied Capelli identities associated with the group
determinant in relation to identities for reducible modules of matrix algebras.
A Capelli identity of irreducible representation type was given in
\cite[Section~3, Theorem~10]{MR3841213} using the subalgebra of the Weyl
algebra associated with each irreducible representation.
After replacing the representation in that paper by its complex conjugate,
we can write the identity in our notation as
\[
(\det X_\varphi)(\det\partial_\varphi)=C_\varphi(0).
\]
Combining this identity with the commutation relations between multiplication
by the determinant and the differential operators yields, for every
$k\in\mathbb Z_{\ge0}$,
\[
(\det X_\varphi)^k(\det\partial_\varphi)^k
=C_\varphi^{(k)}(0).
\]
The first identity in Theorem~\ref{thm:1.9} is obtained by setting
$k=k_\varphi$ in this iterated formula and taking the product over
$\varphi\in I_\Omega$, using the commutativity of operators corresponding
to distinct irreducible representations.
It thus gives an explicit form of the known Capelli identity of irreducible
representation type for the polynomials $P_\Omega$, which generate the
common kernel when $H$ is a proper subgroup.
When the normal core $N$ is trivial, there is only one orbit, and
$P_\Omega=\Theta_G(\bm x)$ and $k_\varphi=d_\varphi$.
Let $\bm\partial=(\partial/\partial x_g)_{g\in G}$, and let
$\Theta_G(\bm\partial)$ denote the operator obtained by replacing each
variable in the group determinant by the corresponding partial derivative.
Since complex conjugation permutes the equivalence classes of irreducible
representations,
\[
\mathcal D_\Omega
=\prod_{\varphi\in\widehat G}(\det\partial_\varphi)^{d_\varphi}
=\Theta_G(\bm\partial).
\]
In this case, the first identity in Theorem~\ref{thm:1.9} reduces to the
known identity of group determinant type
\[
\Theta_G(\bm x)\Theta_G(\bm\partial)
=\prod_{\varphi\in\widehat G}C_\varphi^{(d_\varphi)}(0).
\]

\paragraph{The Cayley identity and differential formulas for powers of the group determinant.}
The classical Cayley identity~\cite[Theorem~2.1]{CaraccioloSokalSportiello2013}
applies to the determinant of a square matrix with independent variable entries.
With the normalization in Corollary~\ref{cor:2.5}, it takes the form
\[
(\det\partial_\varphi)(\det X_\varphi)^r
=\alpha_\varphi^{d_\varphi}
\prod_{i=0}^{d_\varphi-1}(r+i)(\det X_\varphi)^{r-1}
\]
for every positive integer $r$.
For any $s\in\mathbb Z_{\ge0}$, applying this identity $k_\varphi$ times,
starting with the exponent $r=k_\varphi(s+1)$, gives
\[
\begin{aligned}
& (\det\partial_\varphi)^{k_\varphi}
(\det X_\varphi)^{k_\varphi(s+1)}\\
&\qquad=
\alpha_\varphi^{d_\varphi k_\varphi}
\prod_{a=1}^{k_\varphi}\prod_{i=0}^{d_\varphi-1}
(k_\varphi s+a+i)(\det X_\varphi)^{k_\varphi s}.
\end{aligned}
\]
The variables in the matrix blocks corresponding to distinct irreducible
representations are independent.
Taking the product over $\varphi\in I_\Omega$ therefore gives the second
identity in Theorem~\ref{thm:1.9} and the coefficient polynomial $b_\Omega(s)$.

Furthermore, \cite[Theorems~2 and~3]{YI2026} uses Capelli identities of
irreducible representation type and group determinant type to obtain
differential formulas for powers of the irreducible factors and of the
group determinant.
Replacing the representations in that paper by their complex conjugates
to match our notation, and writing its coefficient polynomial as
\[
b_\varphi(t)\coloneqq
\alpha_\varphi^{d_\varphi}
\prod_{j=1}^{d_\varphi}(t+j),
\]
we can express our coefficient polynomial as
\[
b_\Omega(s)
=\prod_{\varphi\in I_\Omega}
\prod_{a=0}^{k_\varphi-1}b_\varphi(k_\varphi s+a).
\]
When $H=G$, each $I_\Omega$ consists of a single irreducible representation,
and $k_\varphi=1$. The second identity in Theorem~\ref{thm:1.9} then
coincides with the differential formula in Theorem~2 of that paper.
When the normal core $N$ is trivial, we have
$P_\Omega=\Theta_G(\bm x)$ and $k_\varphi=d_\varphi$, and the identity
coincides with the differential formula in its Theorem~3.
We describe these differential formulas for the polynomials $P_\Omega$
associated with arbitrary normal cores. For proper subgroups, they are
formulas for the generators of the common kernel.
\section{Proof of Theorem~\ref{thm:1.1}}
In this section, we prove Theorem~\ref{thm:1.1}.
In Subsections~\ref{sec:2.1}, \ref{sec:2.4}, and \ref{sec:2.7},
we assume that $H$ is a proper subgroup of $G$.

\subsection{Preparation~1: Reduction to the normal core}\label{sec:2.1}
Let $G$ be a finite group, let $H$ be a proper subgroup of $G$,
and let $N\coloneqq\Core_G(H)$ be the normal core of $H$ in $G$.
In this subsection, we prove Lemma~\ref{lem:2.1} below.

\begin{lem}\label{lem:2.1}
The invariant rings $K_H^\Delta$ and $K_N^\Delta$ satisfy
\[
 K_H^{\Delta} = K_{N}^{\Delta}.
\]
\end{lem}

The commutator of the polarization operators is defined by
\[
[\Delta_{a}, \Delta_{b}] \coloneqq \Delta_{a} \Delta_{b} - \Delta_{b} \Delta_{a}.
\]
To prove Lemma~\ref{lem:2.1}, we first establish
Lemmas~\ref{lem:2.2} and~\ref{lem:2.3}.

\begin{lem}\label{lem:2.2}
For all $a,b\in G$, we have
\[
[\Delta_{a}, \Delta_{b}] = \Delta_{a b} - \Delta_{b a}.
\]
\end{lem}
\begin{proof}
For every $g\in G$, the identity $\Delta_a(x_g)=x_{ga^{-1}}$ gives
\begin{align*}
\Delta_{a} \Delta_{b}(x_{g})
&= \Delta_{a}(x_{g b^{-1}})
= x_{g b^{-1} a^{-1}}
= x_{g (a b)^{-1}}
= \Delta_{a b}(x_{g}), \\
\Delta_{b} \Delta_{a}(x_{g})
&= \Delta_{b}(x_{g a^{-1}})
= x_{g a^{-1} b^{-1}}
= x_{g (b a)^{-1}}
= \Delta_{b a}(x_{g}).
\end{align*}
Thus, for every $g\in G$, we obtain
\[
 [\Delta_a, \Delta_b](x_g)
 =
 (\Delta_{a b} - \Delta_{b a})(x_g).
\]
The commutator $[\Delta_a,\Delta_b]$ is a derivation because
it is the commutator of two derivations, and
$\Delta_{ab}-\Delta_{ba}$ is also a derivation.
Since they agree on the generators $x_g$ of the polynomial ring
$\mathbb C[\bm x]$, it follows that
\[
 [\Delta_a, \Delta_b]
 =
 \Delta_{a b} - \Delta_{b a}.
\]
\end{proof}

Elements of an invariant ring will be called invariants.

\begin{lem}\label{lem:2.3}
For every invariant $f\in K_H^\Delta$ and all $a,b\in G$, we have
\begin{align*}
 \Delta_{a} f = \Delta_{b^{-1} a b} f.
\end{align*}
\end{lem}
\begin{proof}
First suppose that $b\notin H$.
If $b^{-1}a\notin H$, then
\[
\Delta_{b} f = \Delta_{b^{-1}a} f = 0.
\]
By Lemma~\ref{lem:2.2}, we have
\[
[\Delta_b,\Delta_{b^{-1}a}]
=
\Delta_a-\Delta_{b^{-1}ab},
\]
and hence
\[
 \Delta_{a} f = \Delta_{b^{-1}ab}f.
\]
If $b^{-1}a\in H$, then
\[
a = b(b^{-1}a) \notin H, \quad b^{-1}ab=(b^{-1}a)b \notin H,
\]
so that
\[
 \Delta_af=\Delta_{b^{-1}ab}f=0.
\]
Next suppose that $b\in H$.
Choose an element $c\in G\setminus H$.
Then
\[
b = c^{-1}(cb)
\]
and
\[
 c^{-1} \notin H, \quad c b \notin H.
\]
Applying the case already proved successively to $c^{-1}$ and $cb$, we obtain
\begin{align*}
 \Delta_a f
 =\Delta_{c a c^{-1}} f
 =\Delta_{(cb)^{-1}(c a c^{-1})(cb)}f
 =\Delta_{b^{-1}ab}f.
\end{align*}
\end{proof}

We now prove Lemma~\ref{lem:2.1}.

\begin{proof}[Proof of Lemma~\ref{lem:2.1}]
The inclusion $N\subset H$ immediately gives
$K_N^\Delta\subset K_H^\Delta$.
We prove the reverse inclusion.
Let $f\in K_H^\Delta$.
By the definition of $N$, for every $g\in G\setminus N$,
there exists $t\in G$ such that $g\notin tHt^{-1}$, or equivalently,
\[
t^{-1} g t \notin H.
\]
It follows from Lemma~\ref{lem:2.3} that
\[
 \Delta_g f
 =\Delta_{t^{-1} g t}f
 =0.
\]
Therefore, $f\in K_N^\Delta$.
\end{proof}

\subsection{Preparation~2: Matrix coordinates and the Capelli identity of irreducible representation type}\label{sec:2.2}

Let $\widehat G$ denote a complete set of representatives of the equivalence classes
of irreducible unitary representations of a finite group $G$ over $\mathbb C$.
For each $\varphi\in\widehat G$, set
\[
X_{\varphi} \coloneqq \sum_{g \in G} \varphi(g) x_{g}, \quad
\partial_{\varphi} \coloneqq \sum_{g \in G} \overline{\varphi(g)} \frac{\partial}{\partial x_{g}}, \quad
d_{\varphi} \coloneqq \deg{\varphi}, \quad
\alpha_{\varphi} \coloneqq \frac{|G|}{d_{\varphi}}.
\]
Here, $\overline{\varphi(g)}$ denotes the complex conjugate of the matrix $\varphi(g)$.
Schur's orthogonality relations and the identity
\[
\sum_{\varphi \in \widehat{G}} d_{\varphi}^2 = |G|
\]
show that
\[
\left\{
(X_{\varphi})_{i j}
\mid
\varphi \in \widehat{G}, \
1 \leq i, j \leq d_{\varphi}
\right\}
\]
is a basis of the homogeneous component $\mathbb C[\bm x]_1$ of degree one.
These elements therefore form a new system of linear coordinates on
$\mathbb C[\bm x]$:
\[
\mathbb{C}[\bm{x}]
=
\mathbb{C}[(X_{\varphi})_{i j}
\mid
\varphi \in \widehat{G}, \ 
1 \leq i, j \leq d_{\varphi}].
\]
In this subsection, we record, in our notation, the relations in the Weyl algebra
from \cite[Section~3]{MR3841213} and the Capelli identity of irreducible representation type
from \cite[Theorem~10]{MR3841213}.
Our matrices $X_\varphi$ and $\partial_\varphi$ correspond to those defined
in \cite{MR3841213} for the conjugate representation $\overline\varphi$.

The symbol $\delta_{ij}$ denotes the Kronecker delta.

\begin{lem}\label{lem:2.4}
For all $\varphi,\psi\in\widehat G$, we have
\begin{align*}
[(\partial_\varphi)_{ij},(X_\psi)_{kl}]
=
\begin{cases}
\alpha_\varphi\delta_{ik}\delta_{jl},
&\varphi=\psi,\\
0, &\varphi\ne\psi.
\end{cases}
\end{align*}
\end{lem}
\begin{proof}
Schur's orthogonality relations give
\[
 [(\partial_\varphi)_{ij},(X_\psi)_{kl}]
 =
 \sum_{g\in G}
 \overline{\varphi(g)_{ij}}\psi(g)_{kl}
 =
 \begin{cases}
  \alpha_\varphi\delta_{ik}\delta_{jl},
    & \varphi=\psi,\\
  0,
    & \varphi\ne\psi.
 \end{cases}
\]
\end{proof}

The case $\varphi=\psi$ of Lemma~\ref{lem:2.4}
is given in \cite[Section~3]{MR3841213}.
Corollary~\ref{cor:2.5} below follows from Lemma~\ref{lem:2.4}.

\begin{cor}\label{cor:2.5}
In the coordinate system
$\{(X_\varphi)_{ij}\mid\varphi\in\widehat G,\ 1\leq i,j\leq d_\varphi\}$,
we have
\[
 (\partial_\varphi)_{ij}
 =
 \alpha_\varphi
 \frac{\partial}{\partial (X_\varphi)_{ij}}.
\]
\end{cor}

We next state the Capelli identity of irreducible representation type.
For a $d\times d$ matrix $A=(a_{ij})_{1\leq i,j\leq d}$ with noncommutative entries,
the column determinant is defined by
\[
\Cdet{A} \coloneqq \sum_{\sigma \in S_{d}} \sgn(\sigma) a_{\sigma(1) 1} a_{\sigma(2) 2} \cdots a_{\sigma(d) d}.
\]
For each $\varphi\in\widehat G$, also set
\[
\natural_{\varphi} \coloneqq \diag(d_{\varphi} - 1, d_{\varphi} - 2, \ldots, 0).
\]

\begin{thm}[Capelli identity of irreducible representation type]\label{thm:2.6}
For every $\varphi\in\widehat G$, we have
\[
\Cdet{\left( {}^t X_{\varphi} \partial_{\varphi} + \alpha_{\varphi} \natural_{\varphi} \right)} = \det{X_{\varphi}} \det{\partial_{\varphi}}.
\]
\end{thm}

\subsection{Preparation~3: Description of the common eigenspaces}\label{sec:2.3}

For a function $\lambda\colon G\to\mathbb C$, set
\[
E_{\lambda}
\coloneqq
\left\{
f \in \mathbb{C}[\bm{x}]
\ \middle|\
\Delta_{g} f = \lambda(g) f \ (g \in G)
\right\}.
\]
The character of a representation $\varphi$ is denoted by $\chi_\varphi$.
In this subsection, we prove Lemmas~\ref{lem:2.7} and~\ref{lem:2.8} below.

\begin{lem}\label{lem:2.7}
The common eigenspace $E_\lambda$ is nonzero if and only if
there exist nonnegative integers $m_\varphi$ such that
\[
\lambda
=
\sum_{\varphi \in \widehat{G}}
m_{\varphi} \chi_{\varphi}.
\]
\end{lem}

\begin{lem}\label{lem:2.8}
Suppose that the common eigenspace $E_\lambda$ is nonzero, and write
\[
 \lambda
 =
 \sum_{\varphi \in \widehat{G}}
 m_{\varphi} \chi_{\varphi}
 \quad(m_\varphi \in \mathbb{Z}_{\geq 0}), \qquad
 F_{\lambda} \coloneqq \prod_{\varphi \in \widehat{G}}
 \det{X_{\varphi}}^{m_{\varphi}}.
\]
Then
\[
 E_\lambda
 =
 \mathbb{C} F_{\lambda}.
\]
In particular, the common eigenspace $E_\lambda$ is one-dimensional.
\end{lem}

To prove Lemma~\ref{lem:2.7}, we first state
Lemma~\ref{lem:2.9} and Theorem~\ref{thm:2.10} (Frobenius' formula).

\begin{lem}\label{lem:2.9}
For every $\varphi\in\widehat G$, we have
\begin{align*}
{}^t X_{\varphi} \partial_{\varphi}
=
\sum_{g \in G} \overline{\varphi(g)} \Delta_{g}.
\end{align*}
\end{lem}
\begin{proof}
Since the irreducible representation $\varphi$ is unitary, we have
\begin{align*}
({}^tX_\varphi\partial_\varphi)_{ij}
&=
\sum_{g \in G}
\sum_{h \in G}
\overline{\varphi(g^{-1} h)_{i j}} \ x_{g} \frac{\partial}{\partial x_{h}}\\
&=
\sum_{g \in G} \overline{\varphi(g)_{ij}}
\sum_{h \in G} x_{h g^{-1}} \frac{\partial}{\partial x_{h}}
=
\sum_{g \in G} \overline{\varphi(g)_{ij}} \Delta_{g}.
\end{align*}
\end{proof}

\begin{thm}[Frobenius' formula~{\cite[p.~50]{Frobenius1968gruppen}}, {\cite{Yamaguchi13022026}}]\label{thm:2.10}
For every finite-dimensional representation $\varphi$ of a finite group $G$
over $\mathbb C$ and every $g\in G$, we have
\[
\begin{aligned}
D_{g} \det{\left( \sum_{h \in G} \varphi(h) x_{h} \right)}
&= \Delta_{g} \det{\left( \sum_{h \in G} \varphi(h) x_{h} \right)}\\
&= \chi_{\varphi}(g) \det{\left( \sum_{h \in G} \varphi(h) x_{h} \right)}.
\end{aligned}
\]
\end{thm}

We now prove Lemma~\ref{lem:2.7}.

\begin{proof}[Proof of Lemma~\ref{lem:2.7}]
Suppose that the common eigenspace $E_\lambda$ is nonzero,
and choose $f\in E_\lambda\setminus\{0\}$.
Since $f\in E_\lambda$, Lemma~\ref{lem:2.2} gives, for all $g,h\in G$,
\[
0 = [\Delta_{g}, \Delta_{h}] f = (\Delta_{g h} - \Delta_{h g}) f = \{ \lambda(g h) - \lambda(h g) \} f.
\]
Since $f$ is nonzero, we obtain $\lambda(gh)=\lambda(hg)$,
which implies that $\lambda$ is a class function.
The irreducible characters form a basis of the space of class functions,
so there exist coefficients $m_\psi\in\mathbb C$ such that
\[
 \lambda
 =
 \sum_{\psi \in \widehat{G}}
 m_{\psi} \chi_{\psi}.
\]
Schur's orthogonality relations give, for all $\varphi,\psi\in\widehat G$,
\begin{align*}
\sum_{g \in G}
\overline{\varphi(g)_{ij}} \chi_{\psi}(g)
&=
\sum_{r=1}^{d_\psi}
\sum_{g \in G}
\overline{\varphi(g)_{ij}} \psi(g)_{rr}\\
&=
\begin{cases}
\alpha_{\varphi} \delta_{ij}, & \varphi = \psi, \\
0, & \varphi \neq \psi.
\end{cases}
\end{align*}
Consequently, Lemma~\ref{lem:2.9} yields
\begin{align*}
 ({}^tX_\varphi\partial_\varphi)_{ij}f
 &=
 \sum_{g\in G}
 \overline{\varphi(g)_{ij}}\Delta_gf\\
 &=
 \sum_{g\in G}
 \overline{\varphi(g)_{ij}}\lambda(g)f\\
 &=
 \sum_{\psi \in \widehat{G}}m_{\psi}
 \left(
 \sum_{g \in G}
 \overline{\varphi(g)_{ij}} \chi_{\psi}(g)
 \right)f\\
 &=
 \alpha_{\varphi} m_{\varphi} \delta_{i j} f.
\end{align*}
For the diagonal entries, Corollary~\ref{cor:2.5} gives
\[
 m_\varphi f
 =
 \frac{1}{\alpha_{\varphi}} ({}^{t} X_{\varphi} \partial_{\varphi})_{i i} f
 =
 \sum_{k=1}^{d_\varphi}
 (X_\varphi)_{ki}
 \frac{\partial f}
 {\partial(X_\varphi)_{ki}}.
\]
The operator on the right is the Euler operator in the variables
in the $i$th column of $X_\varphi$.
It multiplies each monomial by its degree in those variables.
Since $f$ is nonzero, its eigenvalue $m_\varphi$ must be a nonnegative integer.
Conversely, suppose that there exist nonnegative integers $m_\varphi$ such that
\[
 \lambda
 =
 \sum_{\varphi \in \widehat{G}}
 m_{\varphi} \chi_{\varphi}.
\]
For the nonzero polynomial
\[
F \coloneqq \prod_{\varphi \in \widehat{G}} \det{X_{\varphi}}^{m_{\varphi}},
\]
Frobenius' formula and the Leibniz rule give
\[
 \Delta_g F
 =
 \lambda(g) F.
\]
Thus,
\[
F \in E_\lambda,
\]
and hence $E_\lambda\ne\{0\}$.
\end{proof}

We now prove Lemma~\ref{lem:2.8}.

\begin{proof}[Proof of Lemma~\ref{lem:2.8}]
The proof of Lemma~\ref{lem:2.7} shows that
\[
 E_\lambda
 \supset
 \mathbb{C} F_{\lambda}.
\]
Choose a nonzero polynomial $f\in E_\lambda$.
The proof of Lemma~\ref{lem:2.7} also gives
\begin{align*}
 ({}^tX_\varphi\partial_\varphi)_{ij}f
 =
 \alpha_\varphi
 m_\varphi\delta_{ij}f.
\end{align*}
We apply the Capelli identity of irreducible representation type,
\[
\det{X_{\varphi}} \det{\partial_{\varphi}} = \Cdet{\left( {}^t X_{\varphi} \partial_{\varphi} + \alpha_{\varphi} \natural_{\varphi} \right)},
\]
to $f$.
Each matrix entry preserves the one-dimensional space $\mathbb Cf$,
so the definition of the column determinant gives
\begin{align*}
\det{X_{\varphi}} \det{\partial_{\varphi}} f
=
\alpha_\varphi^{d_\varphi}
\prod_{j=0}^{d_\varphi-1}
(m_\varphi+j)\,f.
\end{align*}
If $m_\varphi>0$, the scalar on the right is nonzero,
so $f$ is divisible by $\det X_\varphi$.
We may therefore write
\[
 f = f_{1} \det{X_{\varphi}}.
\]
Frobenius' formula and the Leibniz rule give, for every $g\in G$,
\[
\begin{aligned}
\lambda(g) f_{1} \det{X_{\varphi}}
&= \Delta_{g}(f_{1} \det{X_{\varphi}})\\
&= \det{X_{\varphi}} \Delta_{g} f_{1} + \chi_{\varphi}(g) f_{1} \det{X_{\varphi}}.
\end{aligned}
\]
It follows that
\[
\Delta_{g} f_{1} = (\lambda(g) - \chi_{\varphi}(g)) f_{1}.
\]
In other words,
\[
 f_1 \in E_{\lambda-\chi_\varphi}.
\]
Repeating this procedure $m_\varphi$ times for each $\varphi\in\widehat G$, we obtain
\[
f
=
f_{0} F_{\lambda},
\]
where the eigenvalue function of $f_0$ is
\[
\lambda - \sum_{\varphi \in \widehat{G}} m_{\varphi} \chi_{\varphi} = 0.
\]
Thus, for every $g\in G$, we have
\[
 \Delta_{g} f_{0} = 0,
\]
and in particular $\Delta_e f_0=0$.
Since $\Delta_e$ is the Euler operator measuring total degree,
we have $f_0\in\mathbb C$.
Therefore,
\[
 E_\lambda
 \subset
 \mathbb{C} F_{\lambda}.
\]
\end{proof}

\subsection{Preparation~4: Simultaneous diagonalization and decomposition into common eigenspaces}\label{sec:2.4}

For a nonnegative integer $m$, set
\[
K_{H, m}^{\Delta} \coloneqq K_{H}^{\Delta} \cap \mathbb{C}[\bm{x}]_{m},
\]
where $\mathbb C[\bm x]_m$ denotes the homogeneous component
of degree $m$ of $\mathbb C[\bm x]$.
Also set
\[
\Lambda_{H, m} \coloneqq \{ \lambda \colon G \to \mathbb{C} \mid K_{H, m}^{\Delta} \cap E_{\lambda} \neq \{ 0 \} \}.
\]
In this subsection, we prove Lemma~\ref{lem:2.11} below.

\begin{lem}\label{lem:2.11}
We have
\[
K_{H, m}^{\Delta} = \bigoplus_{\lambda \in \Lambda_{H, m}} E_{\lambda}.
\]
\end{lem}

To prove Lemma~\ref{lem:2.11}, we first establish
Lemmas~\ref{lem:2.12} and~\ref{lem:2.13}.

\begin{lem}\label{lem:2.12}
For all $a,b\in G$ and $f\in K_H^\Delta$, we have
\[
 [\Delta_a,\Delta_b]f=0.
\]
\end{lem}
\begin{proof}
Lemmas~\ref{lem:2.2} and~\ref{lem:2.3} give
\[
[\Delta_{a}, \Delta_{b}] f
= \Delta_{a b} f - \Delta_{b a} f
= \Delta_{a b} f - \Delta_{a^{-1} (a b) a} f
= 0.
\]
\end{proof}

\begin{lem}\label{lem:2.13}
For every $a\in G$, we have
\[
 \Delta_a(K_H^\Delta)\subset K_H^\Delta.
\]
\end{lem}
\begin{proof}
For every $g\in G\setminus H$ and $f\in K_H^\Delta$, we have $\Delta_gf=0$.
By Lemma~\ref{lem:2.12}, we also have
\[
 [\Delta_g, \Delta_a]f = 0.
\]
Consequently,
\[
 \begin{aligned}
 \Delta_g(\Delta_af)
 =
 [\Delta_g,\Delta_a]f+\Delta_a(\Delta_gf)
 =0.
 \end{aligned}
\]
Therefore,
\[
 \Delta_af\in
 \bigcap_{g\in G\setminus H}\Ker\Delta_g
 = K_H^\Delta.
\]
\end{proof}

Since the polarization operators $\Delta_g$ preserve degree,
Lemma~\ref{lem:2.13} gives the following corollary.

\begin{cor}\label{cor:2.14}
For every $a\in G$, we have
\[
\Delta_a(K_{H, m}^\Delta)\subset K_{H, m}^\Delta.
\]
\end{cor}

We now prove Lemma~\ref{lem:2.11}.

\begin{proof}[Proof of Lemma~\ref{lem:2.11}]
We first show that, for every $g\in G$,
the operator $\Delta_g$ is diagonalizable on $\mathbb C[\bm x]_m$.
Let $r$ be the order of $g$, and set
\[
\zeta \coloneqq \exp\left(\frac{2\pi \sqrt{-1}}{r}\right).
\]
Let
\[
G = \bigsqcup_{t \in T}t \langle g \rangle
\]
be the decomposition of $G$ into left cosets of the cyclic subgroup $\langle g\rangle$.
For each $t\in T$ and $0\leq k\leq r-1$, set
\[
 z_{t,k}
 :=
 \sum_{j=0}^{r-1}\zeta^{jk}x_{tg^j}.
\]
Since $\Delta_g(x_h)=x_{hg^{-1}}$, we have
\begin{align*}
\Delta_g(z_{t,k})
=
\sum_{j=0}^{r-1}\zeta^{jk}x_{tg^{j-1}}
= \zeta^k \sum_{j=0}^{r-1}\zeta^{jk}x_{tg^j}
= \zeta^k z_{t,k}.
\end{align*}
Thus, each $z_{t,k}$ is an eigenvector of $\Delta_g$.
For each $t\in T$, the transformation
\[
(x_t, x_{tg}, \ldots, x_{tg^{r-1}})
\mapsto
(z_{t, 0}, z_{t, 1}, \ldots, z_{t, r-1})
\]
is a discrete Fourier transform and hence is invertible.
It follows that
\[
 \{z_{t,k}\mid t\in T,\ 0\leq k\leq r-1\}
\]
is a basis of $\mathbb C[\bm x]_1$.
Since $\Delta_g$ is a derivation, for every monomial
\[
M =
\prod_{t\in T}\prod_{k=0}^{r-1}
z_{t,k}^{a_{t,k}} \in \mathbb{C}[\bm{x}]_{m},
\]
we have
\[
 \Delta_g(M)
 =
 \left(
 \sum_{t\in T}\sum_{k=0}^{r-1}
 a_{t,k}\zeta^k
 \right)M.
\]
These monomials form a basis of $\mathbb C[\bm x]_m$.
Therefore, the operator $\Delta_g$ is diagonalizable on $\mathbb C[\bm x]_m$.
By Corollary~\ref{cor:2.14}, the vector space $K_{H,m}^\Delta$
is invariant under $\Delta_g$.
The restriction of a diagonalizable operator to an invariant subspace
is also diagonalizable, so $\Delta_g$ is diagonalizable on $K_{H,m}^\Delta$ as well.
Moreover, Lemma~\ref{lem:2.12} shows that the restrictions of the polarization operators
to $K_{H,m}^\Delta$ commute with one another.
Hence the family $\{\Delta_g\mid g\in G\}$ is simultaneously diagonalizable on this space.
It follows that $K_{H,m}^\Delta$ has a basis of common eigenvectors
of $\{\Delta_g\mid g\in G\}$.
For each vector $f$ in this basis, define its eigenvalue function
$\lambda\colon G\to\mathbb C$ by
\[
\Delta_g f=\lambda(g)f\qquad(g\in G).
\]
The basis vectors corresponding to the same eigenvalue function $\lambda$
span $K_{H,m}^\Delta\cap E_\lambda$.
By the definition of $\Lambda_{H,m}$, we therefore obtain
\[
K_{H,m}^{\Delta}
=
\bigoplus_{\lambda\in\Lambda_{H,m}}
\bigl(K_{H,m}^{\Delta}\cap E_\lambda\bigr).
\]
If $\lambda\in\Lambda_{H,m}$, then
$K_{H,m}^\Delta\cap E_\lambda\ne\{0\}$.
By Lemma~\ref{lem:2.8}, the common eigenspace $E_\lambda$ is one-dimensional.
Since $K_{H,m}^\Delta\cap E_\lambda$ is a nonzero subspace of $E_\lambda$, it follows that
\[
K_{H,m}^{\Delta}\cap E_\lambda=E_\lambda.
\]
Consequently,
\[
K_{H,m}^{\Delta}
=
\bigoplus_{\lambda\in\Lambda_{H,m}}E_\lambda.
\]
\end{proof}
\subsection{Preparation~5: Central idempotents and conjugation orbits}\label{sec:2.5}
In this subsection and the next, let $N$ be an arbitrary normal subgroup of a finite group $G$.
The group algebra of $G$ is
\[
\mathbb{C}G \coloneqq \left\{ \sum_{g \in G} a_{g} g \mid a_{g} \in \mathbb{C} \right\}.
\]
We write $Z(\mathbb{C}G)$ for its center.
We use the same notation for a representation and its $\mathbb C$-linear extension to the group algebra.
For each $\sigma \in \widehat{N}$, put
\[
q_{\sigma} \coloneqq \dfrac{d_{\sigma}}{|N|} \sum_{n \in N} \chi_{\sigma}(n^{-1}) n \in \mathbb{C}N, \qquad d_{\sigma} \coloneqq \deg{\sigma}.
\]
The elements $q_{\sigma}$ are the primitive central idempotents of $\mathbb{C}N$,
and the set $\{ q_{\sigma} \mid \sigma \in \widehat{N} \}$ is a basis of $Z(\mathbb{C}N)$.
For each orbit $\Omega \in \mathcal{O}_{N}$, define
\[
E_{\Omega} := \sum_{\sigma \in \Omega} q_{\sigma} \in Z(\mathbb{C}N).
\]
The aim of this subsection is to prove Lemma~\ref{lem:2.15} below.

\begin{lem}\label{lem:2.15}
For a normal subgroup $N$ of a finite group $G$, we have
\[
Z(\mathbb{C}G) \cap \mathbb{C}N = \bigoplus_{\Omega \in \mathcal{O}_{N}} \mathbb{C} E_{\Omega}.
\]
\end{lem}

To prove Lemma~\ref{lem:2.15}, we first establish Lemmas~\ref{lem:2.16} and \ref{lem:2.17}.

\begin{lem}\label{lem:2.16}
Let $N$ be a normal subgroup of a finite group $G$.
For any $\sigma, \tau \in \widehat{N}$, we have
\begin{align*}
q_{\sigma} q_{\tau} =
\begin{cases}
q_{\sigma}, & \sigma = \tau, \\
0, & \sigma \neq \tau.
\end{cases}
\end{align*}
Moreover,
\[
\sum_{\sigma \in \widehat{N}} q_{\sigma} = e.
\]
\end{lem}
\begin{proof}
A direct calculation gives
\[
q_{\sigma} q_{\tau} = \frac{d_{\sigma} d_{\tau}}{|N|^{2}} \sum_{u \in N} \left( \sum_{a \in N} \chi_{\sigma}(a^{-1}) \chi_{\tau}(u^{-1} a) \right) u.
\]
Schur's orthogonality relations yield
\begin{align*}
\sum_{a \in N} \chi_{\sigma}(a^{-1}) \chi_{\tau}(u^{-1} a) =
\begin{cases}
\dfrac{|N|}{d_{\sigma}} \chi_{\sigma}(u^{-1}), & \sigma = \tau, \\
0, & \sigma \neq \tau,
\end{cases}
\end{align*}
which proves the first assertion.
The character orthogonality relation
\begin{align*}
\sum_{\sigma \in \widehat{N}} d_{\sigma} \chi_{\sigma}(n^{-1})
=
\begin{cases}
|N|, & n = e, \\
0, & n \neq e
\end{cases}
\end{align*}
also gives
\[
\sum_{\sigma \in \widehat{N}} q_{\sigma}
= \frac{1}{|N|} \sum_{n \in N} \left( \sum_{\sigma \in \widehat{N}} d_{\sigma} \chi_{\sigma}(n^{-1}) \right) n
= e.
\]
\end{proof}

\begin{lem}\label{lem:2.17}
The following statements hold:
\begin{enumerate}
\item[(1)] $E_{\Omega} \in Z(\mathbb{C}G) \cap \mathbb{C}N$.
\item[(2)] $E_{\Omega}^{2} = E_{\Omega}$.
\end{enumerate}
\end{lem}
\begin{proof}
We prove (1).
Since $\Omega$ is a $G$-orbit, the map
\[
\Omega \ni \sigma \mapsto g \cdot \sigma \in \Omega
\]
is a bijection for every $g \in G$.
Thus,
\[
g E_{\Omega} g^{-1} = \sum_{\sigma \in \Omega} g q_{\sigma} g^{-1} = \sum_{\sigma \in \Omega} q_{g \cdot \sigma} = E_{\Omega}.
\]
Assertion (2) follows from Lemma~\ref{lem:2.16}.
\end{proof}

\begin{proof}[Proof of Lemma~\ref{lem:2.15}]
By Lemma~\ref{lem:2.17}~(1), the right-hand side is contained in the left-hand side.
We prove the reverse inclusion.
For any $z \in Z(\mathbb{C}G) \cap \mathbb{C}N$, we have $z \in Z(\mathbb{C}N)$.
Hence there are unique coefficients $c_{\sigma} \in \mathbb{C}$ such that
\[
z = \sum_{\sigma \in \widehat{N}} c_{\sigma} q_{\sigma}.
\]
Since $z$ commutes with every element of $G$, the coefficients $c_{\sigma}$ are constant on each $G$-orbit.
Writing $c_{\Omega}$ for the common coefficient on $\Omega$, we obtain
\[
z = \sum_{\Omega \in \mathcal{O}_{N}} c_{\Omega} E_{\Omega}.
\]
Finally, the linear independence of $\{ q_{\sigma} \}_{\sigma \in \widehat{N}}$
implies that of $\{ E_{\Omega} \}_{\Omega \in \mathcal{O}_{N}}$.
\end{proof}

\subsection{Preparation~6: Character vanishing conditions and the description of exponents}\label{sec:2.6}

For each orbit $\Omega \in \mathcal{O}_{N}$, put
\[
I_{\Omega} \coloneqq \{ \varphi \in \widehat{G} \mid \varphi|_{N} \text{ has an irreducible constituent in } \Omega \}.
\]
Also, let
\[
\delta_{\Omega} \coloneqq \gcd{\{ d_{\varphi} \mid \varphi \in I_{\Omega} \}}.
\]
The aim of this subsection is to prove Lemma~\ref{lem:2.18} below.

\begin{lem}\label{lem:2.18}
Let $N$ be a normal subgroup of a finite group $G$, and let
$m_\varphi\in\mathbb Z_{\geq0}$ for each $\varphi\in\widehat G$.
The condition
\[
\sum_{\varphi \in \widehat{G}} m_{\varphi} \chi_{\varphi}(g) = 0
\]
holds for every $g \in G \setminus N$ if and only if,
for each $\Omega \in \mathcal{O}_{N}$, there exists $k_{\Omega} \in \mathbb{Z}_{\geq 0}$ such that
\[
m_{\varphi} = k_{\Omega} \dfrac{d_{\varphi}}{\delta_{\Omega}}
\]
for every $\varphi \in I_{\Omega}$.
\end{lem}

For each irreducible representation $\varphi \in \widehat{G}$, put
\[
p_{\varphi} \coloneqq \dfrac{d_{\varphi}}{|G|} \sum_{g \in G} \chi_{\varphi}(g^{-1}) g \in \mathbb{C}G.
\]
To prove Lemma~\ref{lem:2.18}, we first establish Lemmas~\ref{lem:2.19}--\ref{lem:2.21}.

\begin{lem}\label{lem:2.19}
For each orbit $\Omega \in \mathcal{O}_{N}$, we have
\[
E_{\Omega} = \sum_{\varphi \in I_{\Omega}} p_{\varphi}.
\]
\end{lem}
\begin{proof}
By Lemma~\ref{lem:2.17}~(1), we have $E_\Omega\in Z(\mathbb C G)$.
Since $\{p_\varphi\mid\varphi\in\widehat G\}$ is a basis of $Z(\mathbb{C}G)$,
there are unique coefficients $a_{\Omega, \varphi} \in \mathbb{C}$ such that
\[
E_{\Omega} = \sum_{\varphi \in \widehat{G}} a_{\Omega, \varphi} p_{\varphi}.
\]
Applying Lemma~\ref{lem:2.16} with $N = G$, we see that the elements $p_\varphi$ are pairwise orthogonal idempotents.
Together with Lemma~\ref{lem:2.17}~(2), this gives
\[
\sum_{\varphi \in \widehat{G}} a_{\Omega, \varphi} p_{\varphi} = E_{\Omega} = E_{\Omega}^{2} = \sum_{\varphi \in \widehat{G}} a_{\Omega, \varphi}^{2} p_{\varphi}.
\]
Hence $a_{\Omega, \varphi} \in \{ 0, 1 \}$.
For $\varphi \in \widehat{G}$ and $\sigma \in \widehat{N}$, put
\[
m_{\varphi, \sigma} := \dfrac{1}{|N|} \sum_{n \in N} \chi_{\varphi}(n) \chi_{\sigma}(n^{-1}).
\]
This is the multiplicity of $\chi_{\sigma}$ in $\chi_{\varphi}|_{N}$.
By the definition of $q_{\sigma}$, we have
\[
\Tr{\varphi(q_{\sigma})}
= \frac{d_{\sigma}}{|N|} \sum_{n \in N} \chi_{\sigma}(n^{-1}) \chi_{\varphi}(n)
= d_{\sigma} m_{\varphi, \sigma},
\]
and therefore
\[
\Tr{\varphi(E_{\Omega})} = \sum_{\sigma \in \Omega} d_{\sigma} m_{\varphi, \sigma}.
\]
The right-hand side is a sum of nonnegative integers, and each $d_\sigma$ is positive.
It follows that
\[
\Tr{\varphi(E_{\Omega})} > 0 \iff \varphi \in I_{\Omega}.
\]
On the other hand, since $p_{\psi} \in Z(\mathbb{C}G)$,
Schur's lemma gives a scalar $c \in \mathbb{C}$ such that
\[
\varphi(p_{\psi}) = c I_{d_{\varphi}}.
\]
The orthogonality relations for irreducible characters yield
\begin{align*}
\Tr{\varphi(p_{\psi})}
= \frac{d_{\psi}}{|G|} \sum_{g \in G} \chi_{\psi}(g^{-1}) \chi_{\varphi}(g)
=
\begin{cases}
d_{\varphi}, & \psi = \varphi, \\
0, & \psi \neq \varphi,
\end{cases}
\end{align*}
so that
\begin{align*}
\varphi(p_{\psi})
=
\begin{cases}
I_{d_{\varphi}}, & \psi = \varphi, \\
0, & \psi \neq \varphi.
\end{cases}
\end{align*}
Consequently,
\[
\Tr{\varphi(E_{\Omega})} = a_{\Omega, \varphi} d_{\varphi}.
\]
Combining these observations, we obtain
\[
a_{\Omega, \varphi} = 1 \iff \varphi \in I_{\Omega},
\]
and hence
\[
E_{\Omega} = \sum_{\varphi \in I_{\Omega}} p_{\varphi}.
\]
\end{proof}

\begin{lem}\label{lem:2.20}
The family $\{I_\Omega \}_{\Omega\in\mathcal{O}_N}$ is a partition of $\widehat G$.
\end{lem}
\begin{proof}
The character orthogonality relations give
\[
\sum_{\varphi \in \widehat{G}} p_{\varphi} = e.
\]
Thus, Lemmas~\ref{lem:2.16} and \ref{lem:2.19} imply
\[
\sum_{\Omega \in \mathcal O_{N}} \sum_{\varphi \in I_{\Omega}} p_{\varphi}
=
\sum_{\Omega \in \mathcal O_{N}} E_{\Omega}
=
\sum_{\Omega \in \mathcal O_{N}} \sum_{\sigma \in \Omega} q_{\sigma}
=
e
=
\sum_{\varphi \in \widehat{G}} p_{\varphi}.
\]
The coefficient of $p_\varphi$ on the left-hand side is the number of sets $I_\Omega$ containing $\varphi$.
Since the set $\{p_\varphi\}_{\varphi\in\widehat G}$ is linearly independent
and the coefficient of each $p_\varphi$ on the right-hand side is $1$,
each $\varphi\in\widehat G$ belongs to exactly one $I_\Omega$.
Moreover, each $I_\Omega$ is nonempty because $E_\Omega\ne 0$.
\end{proof}

For each $\Omega\in\mathcal{O}_N$, put
\[
\eta_\Omega:=\sum_{\varphi\in I_\Omega}d_\varphi\chi_\varphi.
\]

\begin{lem}\label{lem:2.21}
The set $\{\eta_\Omega\}_{\Omega\in\mathcal{O}_N}$ is a basis of the space of class functions on $G$ that vanish on $G \setminus N$.
\end{lem}
\begin{proof}
Let $L(G)$ be the vector space of complex-valued class functions on the finite group $G$.
Consider the linear map
\[
\Phi\colon L(G)\to Z(\mathbb{C}G), \qquad
\Phi(\mu) \coloneqq \frac{1}{|G|}\sum_{g\in G}\mu(g^{-1})g.
\]
If $\mu$ is a class function, then $\mu(g^{-1})$ is constant as $g$ ranges over a conjugacy class,
so $\Phi(\mu) \in Z(\mathbb{C}G)$.
By the definition of $p_\varphi$, we have
\[
\Phi(\chi_\varphi) = \frac{1}{d_\varphi}p_\varphi.
\]
Since $\{\chi_\varphi\}_{\varphi\in\widehat G}$ and $\{p_\varphi\}_{\varphi\in\widehat G}$
are bases of $L(G)$ and $Z(\mathbb{C}G)$, respectively, the map $\Phi$ is a linear isomorphism.
The definition of $\Phi$ gives
\[
\mu(g)=0\quad(g\in G\setminus N)
\quad\Longleftrightarrow\quad
\Phi(\mu)\in Z(\mathbb{C}G)\cap\mathbb{C}N.
\]
Moreover, Lemma~\ref{lem:2.19} yields
\[
\Phi(\eta_\Omega)
=\sum_{\varphi\in I_\Omega}d_\varphi\Phi(\chi_\varphi)
=\sum_{\varphi\in I_\Omega}p_\varphi
=E_\Omega.
\]
By Lemma~\ref{lem:2.15}, the set $\{E_\Omega\}_{\Omega\in\mathcal{O}_N}$
is a basis of $Z(\mathbb{C}G)\cap\mathbb{C}N$.
Since $\Phi$ is a linear isomorphism, the equivalence above shows that
$\{\eta_\Omega\}_{\Omega\in\mathcal{O}_N}$ is a basis of the space of class functions that vanish on $G\setminus N$.
\end{proof}

\begin{proof}[Proof of Lemma~\ref{lem:2.18}]
Suppose that the function
\[
\lambda \coloneqq
\sum_{\varphi\in\widehat G}m_\varphi\chi_\varphi
\]
satisfies $\lambda(g) = 0$ for every $g\in G\setminus N$.
By Lemma~\ref{lem:2.21}, there are coefficients $c_\Omega\in\mathbb{C}$ such that
\[
\lambda
=\sum_{\Omega\in\mathcal{O}_N}c_\Omega\eta_\Omega
=\sum_{\Omega\in\mathcal{O}_N}
 \sum_{\varphi\in I_\Omega}c_\Omega d_\varphi\chi_\varphi.
\]
By Lemma~\ref{lem:2.20}, the family $\{I_\Omega\}_{\Omega\in\mathcal{O}_N}$ is a partition of $\widehat G$.
Comparing coefficients and using the linear independence of the irreducible characters, we obtain
\[
m_\varphi=c_\Omega d_\varphi
\qquad
(\Omega\in\mathcal{O}_N,\ \varphi\in I_\Omega).
\]
Fix $\Omega \in \mathcal{O}_N$.
Since $\delta_\Omega$ is the greatest common divisor of the integers $d_\varphi$ with $\varphi \in I_\Omega$,
B\'ezout's identity gives integers $u_\varphi\in\mathbb{Z}$ such that
\[
\delta_\Omega=\sum_{\varphi\in I_\Omega}u_\varphi d_\varphi.
\]
It follows that
\[
c_\Omega\delta_\Omega
=\sum_{\varphi\in I_\Omega}u_\varphi m_\varphi\in\mathbb{Z}.
\]
Since $I_\Omega \neq \emptyset$, we may choose $\varphi_0\in I_\Omega$, which gives
\[
c_\Omega=\frac{m_{\varphi_0}}{d_{\varphi_0}}\geq0.
\]
Thus, setting $k_\Omega\coloneqq c_\Omega\delta_\Omega$, we have
$k_\Omega\in\mathbb{Z}_{\geq 0}$ and
\[
m_\varphi=k_\Omega\frac{d_\varphi}{\delta_\Omega}
\qquad(\varphi\in I_\Omega).
\]
Conversely, suppose that for each $\Omega\in\mathcal{O}_N$, there exists
$k_\Omega\in\mathbb{Z}_{\geq0}$ such that
$m_\varphi=k_\Omega d_\varphi/\delta_\Omega$ for every $\varphi\in I_\Omega$.
By Lemma~\ref{lem:2.20}, we can write
\[
\lambda
=\sum_{\Omega\in\mathcal{O}_N}
 \sum_{\varphi\in I_\Omega}
 k_\Omega\frac{d_\varphi}{\delta_\Omega}\chi_\varphi
=\sum_{\Omega\in\mathcal{O}_N}
 \frac{k_\Omega}{\delta_\Omega}\eta_\Omega.
\]
By Lemma~\ref{lem:2.21}, each $\eta_\Omega$ vanishes on $G\setminus N$.
Hence $\lambda$ also vanishes on $G\setminus N$.
\end{proof}

\subsection{Proof of Theorem~\ref{thm:1.1}}\label{sec:2.7}
Let $H$ be a proper subgroup of $G$, and again put $N=\Core_G(H)$.
For each $\Omega \in \mathcal{O}_{N}$, define
\[
P_{\Omega} \coloneqq \prod_{\varphi \in I_{\Omega}} \det{X_{\varphi}}^{\frac{d_{\varphi}}{\delta_{\Omega}}}.
\]
We now prove Theorem~\ref{thm:1.1}.

\begin{proof}[Proof of Theorem~\ref{thm:1.1}]
By Lemma~\ref{lem:2.1}, we have $K_H^\Delta=K_N^\Delta$.
We first determine $K_N^\Delta$.
Since each $\Delta_g$ preserves degree, we have
\[
K_N^\Delta=\bigoplus_{m\ge0}K_{N,m}^\Delta.
\]
By Lemmas~\ref{lem:2.7}, \ref{lem:2.8}, and \ref{lem:2.11}, the space $K_N^\Delta$ is spanned by the polynomials
\[
F_{\lambda} = \prod_{\varphi \in \widehat{G}}
\det{X_\varphi}^{m_\varphi},
\qquad
\lambda=\sum_{\varphi\in\widehat G}m_\varphi\chi_\varphi,
\qquad
m_\varphi\in\mathbb Z_{\geq 0}
\]
for which $\lambda(g) = 0$ for every $g\in G \setminus N$.
By Lemma~\ref{lem:2.18}, this condition is equivalent to the existence,
for each $\Omega\in\mathcal O_N$, of $k_\Omega\in\mathbb Z_{\geq 0}$ such that
\[
m_\varphi=k_\Omega\frac{d_\varphi}{\delta_\Omega}
\qquad(\varphi\in I_\Omega).
\]
Thus, these polynomials $F_\lambda$ are precisely those of the form
\[
\prod_{\Omega\in\mathcal O_N}
\left(
  \prod_{\varphi\in I_\Omega}
  \det{X_\varphi}^{\frac{d_\varphi}{\delta_\Omega}}
\right)^{k_\Omega}
=
\prod_{\Omega \in \mathcal{O}_{N}} P_{\Omega}^{k_{\Omega}}.
\]
Consequently,
\[
K_H^\Delta=K_N^\Delta
=\mathbb C[P_\Omega\mid\Omega\in\mathcal O_N].
\]
Next, we prove $K_{H}^{\Delta} = K_{H}^{D} = K_{N}^{D}$.
Frobenius' formula gives
\[
D_g \det{X_\varphi} = \Delta_g \det{X_\varphi}
=\chi_\varphi(g) \det{X_\varphi}
\]
for each $\det{X_\varphi}$.
Since both operators are derivations, the Leibniz rule gives $D_gP_\Omega=\Delta_gP_\Omega$.
Hence
\[
\mathbb C[P_\Omega\mid\Omega\in\mathcal O_N]
\subset K_N^D\subset K_H^D.
\]
On the other hand, consider the ring automorphism $\iota \colon \mathbb{C}[\bm{x}] \to \mathbb{C}[\bm{x}]$
defined by $\iota(x_h)\coloneqq x_{h^{-1}}$.
Then $\iota^2=\mathrm{id}$.
The operators $\iota\circ\Delta_{g^{-1}}\circ\iota$ and $D_g$ are both derivations,
and on the generators we have
\[
(\iota\circ\Delta_{g^{-1}}\circ\iota)(x_h)
=\iota(x_{h^{-1}g})
=x_{g^{-1}h}
=D_g(x_h).
\]
Therefore,
\[
\iota\circ\Delta_{g^{-1}}\circ\iota=D_g.
\]
Since $\iota^2=\mathrm{id}$, it follows that
\[
D_g\circ\iota=\iota\circ\Delta_{g^{-1}}.
\]
Thus, for every $f\in\mathbb C[\bm{x}]$, we have
\[
f\in K_H^\Delta
\quad\Longleftrightarrow\quad
\iota(f)\in K_H^D,
\]
and consequently
\[
\iota(K_H^\Delta)=K_H^D.
\]
Applying $\iota$ to the inclusion $K_H^\Delta \subset K_H^D$ obtained above
and using $\iota^2 = \mathrm{id}$, we get
\[
K_H^D
=\iota(K_H^\Delta)
\subset \iota(K_H^D)
=K_H^\Delta.
\]
Hence $K_H^D=K_H^\Delta$.
Combining these equalities, we obtain
\[
K_H^D=K_H^\Delta=K_N^D=K_N^\Delta
=\mathbb C[P_\Omega\mid\Omega\in\mathcal O_N].
\]
\end{proof}
\section{Proof of Theorem~\ref{thm:1.2}}
In this section, we prove Theorem~\ref{thm:1.2}.

\begin{proof}[Proof of Theorem~\ref{thm:1.2}]
The entries of the matrices $X_\varphi$ $(\varphi\in\widehat G)$,
\[
\left\{
(X_{\varphi})_{i j}
\mid
\varphi \in \widehat{G}, \
1 \leq i, j \leq d_{\varphi}
\right\},
\]
form a system of linear coordinates for $\mathbb C[\bm{x}]$.
Moreover, the family $\{I_\Omega\mid\Omega\in\mathcal O_N\}$ is a partition of
$\widehat G$.
Thus, using independent variables $x_\Omega$ $(\Omega\in\mathcal O_N)$,
we may make the substitution
\[
X_\varphi=\operatorname{diag}(x_\Omega, 1, 1, \ldots, 1)
\qquad(\varphi \in I_\Omega).
\]
Under this substitution,
\[
P_\Omega\longmapsto x_\Omega^{r_\Omega},
\qquad
r_\Omega:=
\sum_{\varphi\in I_\Omega}\frac{d_\varphi}{\delta_\Omega}
\in\mathbb Z_{\geq 1}.
\]
Suppose that a polynomial
$F \in \mathbb{C}[T_\Omega\mid\Omega\in\mathcal{O}_N]$ satisfies
\[
F\bigl((P_\Omega)_{\Omega\in\mathcal O_N}\bigr)=0.
\]
The substitution above gives
\[
F\bigl((x_\Omega^{r_\Omega})_{\Omega\in\mathcal O_N}\bigr)=0.
\]
Since each $r_\Omega$ is a positive integer,
the powers $x_\Omega^{r_\Omega}$ independently take all values in $\mathbb C$.
Hence $F$ vanishes at every tuple of complex numbers, and therefore $F=0$.
Thus $\{P_\Omega \mid \Omega \in \mathcal{O}_N\}$ is algebraically independent
over $\mathbb C$.
\end{proof}

\section{Proof of Theorem~\ref{thm:1.4}}
In this section, we prove Theorem~\ref{thm:1.4}.
We first establish Lemma~\ref{lem:4.1}.

\begin{lem}\label{lem:4.1}
Let $H$ be a proper subgroup of a finite group $G$, and set
$N \coloneqq \Core_G(H)$.
Then
\[
N
=
\{g\in G\mid
  \Delta_g f\ne0 \text{ for some } f\in K_H\}.
\]
\end{lem}
\begin{proof}
For each $\Omega\in\mathcal O_N$, set
\[
\lambda_\Omega
\coloneqq
\sum_{\varphi\in I_\Omega}
\frac{d_\varphi}{\delta_\Omega}\chi_\varphi
=
\frac{\eta_\Omega}{\delta_\Omega}.
\]
Since $\delta_\Omega\geq1$, Lemma~\ref{lem:2.21} shows that
$\{\lambda_\Omega\}_{\Omega\in\mathcal O_N}$ is a basis of the space of
class functions on $G$ that vanish on $G\setminus N$.
By the definition of $P_\Omega$, Theorem~\ref{thm:2.10},
and the Leibniz rule, we also have
\[
\Delta_gP_\Omega
=
\left(
\sum_{\varphi\in I_\Omega}
\frac{d_\varphi}{\delta_\Omega}\chi_\varphi(g)
\right)P_\Omega
=
\lambda_\Omega(g)P_\Omega
\qquad(g\in G).
\]
First, suppose that $g\in G\setminus N$.
By Theorem~\ref{thm:1.1}, we have $K_H=K_N^\Delta$.
The definition of the invariant ring therefore gives
$\Delta_g f=0$ for every $f\in K_H$.
Conversely, suppose that $g\in N$, and let $C_g$ be its conjugacy class in $G$.
Since $N$ is normal in $G$, we have $C_g\subset N$.
Let $\mathbf 1_{C_g}\colon G\to\mathbb C$ be the function that equals $1$
on $C_g$ and $0$ elsewhere.
This is a class function that vanishes on $G\setminus N$, so we may write
\[
\mathbf 1_{C_g}
=
\sum_{\Omega\in\mathcal O_N}c_\Omega\lambda_\Omega
\]
for some $c_\Omega\in\mathbb C$.
Evaluating both sides at $g$ gives
\[
1
=
\sum_{\Omega\in\mathcal O_N}c_\Omega\lambda_\Omega(g).
\]
Thus $\lambda_{\Omega_0}(g)\ne0$ for some $\Omega_0\in\mathcal O_N$.
By Theorem~\ref{thm:1.1}, we have $P_{\Omega_0}\in K_H$.
Moreover, substituting $x_e=1$ and $x_h=0$ $(h\ne e)$ gives
$P_{\Omega_0}=1$, so $P_{\Omega_0}\ne0$.
Consequently,
\[
\Delta_gP_{\Omega_0}
=
\lambda_{\Omega_0}(g)P_{\Omega_0}
\ne0.
\]
This proves the assertion.
\end{proof}

We now prove Theorem~\ref{thm:1.4}.

\begin{proof}[Proof of Theorem~\ref{thm:1.4}]
First, suppose that $N_1\subset N_2$.
Then $G\setminus N_2\subset G\setminus N_1$.
The definition of the invariant rings and Theorem~\ref{thm:1.1} give
\[
K_{H_1}
=
K_{N_1}^\Delta
\subset
K_{N_2}^\Delta
=
K_{H_2}.
\]
Conversely, suppose that $K_{H_1}\subset K_{H_2}$.
For any $g\in N_1$, Lemma~\ref{lem:4.1} gives an element $f\in K_{H_1}$
such that $\Delta_g f\ne0$.
By assumption, we also have $f\in K_{H_2}$, so another application of
Lemma~\ref{lem:4.1} yields $g\in N_2$.
Hence $N_1\subset N_2$.
We have therefore proved
\[
K_{H_1}\subset K_{H_2}
\iff
N_1\subset N_2.
\]
Applying this equivalence to the inclusions in both directions also gives
\[
K_{H_1}=K_{H_2}
\iff
N_1=N_2.
\]
\end{proof}

\section{Proof of Theorem~\ref{thm:1.5}}
In this section, we prove Theorem~\ref{thm:1.5}.

\begin{proof}[Proof of Theorem~\ref{thm:1.5}]
By Lemma~\ref{lem:2.20}, the family
$\{I_\Omega\}_{\Omega\in\mathcal O_N}$ is a partition of $\widehat G$.
Since $k_\varphi\delta_\Omega=d_\varphi$ for $\varphi\in I_\Omega$, we obtain
\[
\begin{aligned}
\prod_{\Omega\in\mathcal O_N}P_\Omega^{\delta_\Omega}
=
\prod_{\Omega\in\mathcal O_N}
\prod_{\varphi\in I_\Omega}
(\det X_\varphi)^{k_\varphi\delta_\Omega}
=
\prod_{\varphi\in\widehat G}(\det X_\varphi)^{d_\varphi}
=\Theta_G(\bm x).
\end{aligned}
\]
The last equality follows from Frobenius' group determinant theorem.
\end{proof}

\section{Proof of Theorem~\ref{thm:1.7}}
In this section, we prove Theorem~\ref{thm:1.7}.

\begin{proof}[Proof of Theorem~\ref{thm:1.7}]
Since $G$ is abelian, we have $N=H$, and the conjugation action of $G$
on $\widehat H$ is trivial.
Thus
\[
\mathcal O_H
=
\{\{\sigma\}\mid\sigma\in\widehat H\}.
\]
Moreover, all irreducible representations of $G$ and $H$ are one-dimensional.
Fix an irreducible character $\sigma\in\widehat H$.
Then
\[
I_{\{\sigma\}}
=
\{\varphi\in\widehat G\mid\varphi|_H=\sigma\}.
\]
By Lemma~\ref{lem:2.20}, the set $I_{\{\sigma\}}$ is nonempty, so we may
choose $\varphi_\sigma\in\widehat G$ such that $\varphi_\sigma|_H=\sigma$.
Since $d_\varphi=1$ for every $\varphi\in I_{\{\sigma\}}$, we also have
$\delta_{\{\sigma\}}=1$.
Let $\pi\colon G\to G/H$ be the natural projection.
An irreducible representation $\varphi\in\widehat G$ satisfies
$\varphi|_H=\sigma$ if and only if $\varphi\varphi_\sigma^{-1}$ is trivial
on $H$.
Every character of $G$ that is trivial on $H$ uniquely determines a character
of $G/H$, so the map
\[
\widehat{G/H}\longrightarrow I_{\{\sigma\}},
\qquad
\psi\longmapsto\varphi_\sigma(\psi\circ\pi)
\]
is a bijection.
The definition of $P_{\{\sigma\}}$ therefore gives
\[
\begin{aligned}
P_{\{\sigma\}}
&=
\prod_{\varphi\in I_{\{\sigma\}}}
\left(\sum_{g\in G}\varphi(g)x_g\right)\\
&=
\prod_{\psi\in\widehat{G/H}}
\left(
\sum_{g\in G}\varphi_\sigma(g)\psi(gH)x_g
\right)\\
&=
\prod_{\psi\in\widehat{G/H}}
\left(
\sum_{tH\in G/H}\psi(tH)
\sum_{g\in tH}\varphi_\sigma(g)x_g
\right)\\
&=
\prod_{\psi\in\widehat{G/H}}
\left(
\sum_{tH\in G/H}\psi(tH)y_{tH}^\sigma
\right)\\
&=
\Theta_{G/H}(\bm y^\sigma).
\end{aligned}
\]
The last equality follows from Frobenius' group determinant theorem
for the finite abelian group $G/H$.
\end{proof}

\section{Proof of Theorem~\ref{thm:1.9}}
In this section, we prove Theorem~\ref{thm:1.9}.

\begin{proof}[Proof of Theorem~\ref{thm:1.9}]
Fix an orbit $\Omega\in\mathcal O_N$.
By Lemma~\ref{lem:2.4}, if $\varphi\ne\psi$, then every entry of
$X_\varphi$ and $\partial_\varphi$ commutes with every entry of
$X_\psi$ and $\partial_\psi$.
We first prove the first identity.
By Theorem~\ref{thm:2.6}, we have
\[
\det X_\varphi\det\partial_\varphi=C_\varphi(0).
\]
Moreover, Lemma~\ref{lem:2.4} and the multilinearity of the determinant give
\[
\begin{aligned}
\bigl[
({}^tX_\varphi\partial_\varphi)_{ij},
\det X_\varphi
\bigr]
&=
\alpha_\varphi\delta_{ij}\det X_\varphi,\\
\bigl[
({}^tX_\varphi\partial_\varphi)_{ij},
\det\partial_\varphi
\bigr]
&=
-\alpha_\varphi\delta_{ij}\det\partial_\varphi.
\end{aligned}
\]
Applying the latter commutation relation to each term of the column determinant,
we obtain
\[
C_\varphi(u)(\det\partial_\varphi)
=
(\det\partial_\varphi)C_\varphi(u-1).
\]
Thus, for $k\in\mathbb Z_{\ge1}$,
\[
\begin{aligned}
(\det X_\varphi)^k(\det\partial_\varphi)^k
&=
(\det X_\varphi)^{k-1}
C_\varphi(0)(\det\partial_\varphi)^{k-1}\\
&=
(\det X_\varphi)^{k-1}
(\det\partial_\varphi)^{k-1}C_\varphi(1-k).
\end{aligned}
\]
Induction on $k$ now yields
\[
(\det X_\varphi)^k(\det\partial_\varphi)^k
=
C_\varphi(0)C_\varphi(-1)\cdots C_\varphi(1-k)
=
C_\varphi^{(k)}(0).
\]
We apply this identity with $k=k_\varphi$.
Since the operators corresponding to distinct irreducible representations
commute with one another, we obtain
\[
\begin{aligned}
P_\Omega\mathcal D_\Omega
&=
\prod_{\varphi\in I_\Omega}
(\det X_\varphi)^{k_\varphi}
(\det\partial_\varphi)^{k_\varphi}\\
&=
\prod_{\varphi\in I_\Omega}
C_\varphi^{(k_\varphi)}(0).
\end{aligned}
\]
This proves the first identity.
We next prove the second identity.
Let $s\in\mathbb Z_{\ge0}$.
For each $\varphi\in I_\Omega$ and $1\leq i,j\leq d_\varphi$, the operator
$({}^tX_\varphi\partial_\varphi)_{ij}$ is a derivation.
The commutation relations above and the Leibniz rule therefore give
\[
({}^tX_\varphi\partial_\varphi)_{ij}
(P_\Omega^{s+1})
=
\alpha_\varphi k_\varphi(s+1)\delta_{ij}P_\Omega^{s+1}.
\]
Thus these operators preserve the one-dimensional space
$\mathbb C P_\Omega^{s+1}$.
On this space, the off-diagonal entries act as zero, while each diagonal entry
acts as multiplication by $\alpha_\varphi k_\varphi(s+1)$.
The definition of the column determinant gives
\[
C_\varphi(u)(P_\Omega^{s+1})
=
\alpha_\varphi^{d_\varphi}
\prod_{i=0}^{d_\varphi-1}
\bigl(k_\varphi(s+1)+u+i\bigr)
P_\Omega^{s+1}.
\]
Applying this formula with $u=0,-1,\ldots,1-k_\varphi$, we obtain
\begin{align*}
C_\varphi^{(k_\varphi)}(0)(P_\Omega^{s+1})
&=
\alpha_\varphi^{d_\varphi k_\varphi}
\prod_{j=0}^{k_\varphi-1}
\prod_{i=0}^{d_\varphi-1}
\bigl(k_\varphi(s+1)-j+i\bigr)
P_\Omega^{s+1}\\
&=
\alpha_\varphi^{d_\varphi k_\varphi}
\prod_{a=1}^{k_\varphi}
\prod_{i=0}^{d_\varphi-1}
(k_\varphi s+a+i)
P_\Omega^{s+1}.
\end{align*}
In the last equality, we have substituted $a=k_\varphi-j$.
Taking the product over $\varphi\in I_\Omega$ and using the first identity,
we obtain
\[
\begin{aligned}
P_\Omega\mathcal D_\Omega(P_\Omega^{s+1})
&=
\left(
\prod_{\varphi\in I_\Omega}
C_\varphi^{(k_\varphi)}(0)
\right)
(P_\Omega^{s+1})\\
&=
b_\Omega(s)P_\Omega^{s+1}.
\end{aligned}
\]
Since $P_\Omega$ is nonzero and $\mathbb C[\bm x]$ is an integral domain,
we may cancel $P_\Omega$ from both sides to obtain
\[
\mathcal D_\Omega(P_\Omega^{s+1})
=
b_\Omega(s)P_\Omega^s.
\]
\end{proof}

\clearpage

\bibliography{reference}

\medskip
\begin{flushleft}
Faculty of Education,
University of Miyazaki,
1-1 Gakuen Kibanadai-nishi,
Miyazaki 889-2192,
Japan \\
{\it Email address}, Naoya Yamaguchi: n-yamaguchi@miyazaki-u.ac.jp \\
\end{flushleft}

\end{document}